\documentclass[a4paper,12pt]{article}%
\scrollmode

\usepackage{amsmath,amssymb,amsfonts}
\usepackage{graphicx}
\newtheorem{theorem}{Theorem}[section]

\newtheorem{lemma}[theorem]{Lemma}
\newtheorem{proposition}[theorem]{Proposition}
\newtheorem{corollary}[theorem]{Corollary}
\newtheorem{remark}[theorem]{Remark}
\usepackage{setspace}
\usepackage{color}
 
\title{Long range trap models on $\mathbb{Z}$ and quasistable processes\thanks{
Paper accepted for publication in the {\it Journal of Theoretical Probability.}}}
\author{W. Barreto-Souza\footnote{Email: wagnerbs85@gmail.com}\,\,\thanks{Supported by a CNPq doctoral fellowship.}\, 
and L.R.G. Fontes\footnote{Email: lrfontes@usp.br}
\thanks{Partially supported by CNPq grant 305760/2010-6, and FAPESP grant 2009/52379-8.}\\\\
Instituto de Matem\'atica e Estat\' \i stica\\ Universidade de S\~ao Paulo}
\date{}

\begin{document}

\maketitle

\begin{abstract}

Let $\mathcal X=\{\mathcal X_t:\, t\geq0,\, \mathcal X_0=0\}$ be a mean zero $\beta$-stable random walk on $\mathbb{Z}$ 
with inhomogeneous jump rates $\{\tau_i^{-1}: i\in\mathbb{Z}\}$, with $\beta\in(1,2]$ and $\{\tau_i: i\in\mathbb{Z}\}$ 
a family of independent random variables with common marginal distribution in the basin of attraction of an 
$\alpha$-stable law, $\alpha\in(0,1)$. In this paper we derive results about the long time behavior of this process, 
in particular its scaling limit, given by a $\beta$-stable process time-changed by 
the inverse of another process, involving the local time of the $\beta$-stable process and an independent $\alpha$-stable 
subordinator; we call the resulting process a {\it quasistable process}. 
Another such result concerns {\em aging}. We obtain an (integrated) 
aging result for $\mathcal X$. \\


\end{abstract}

\section{Introduction}

Trap models have been introduced in the physics literature as simple models of disordered systems where long time memory 
effects like aging and localization can be established and understood on a rigorous basis. 
See for instance~\cite{bouchaud1992},~\cite{bouchaud-dean1995} and~\cite{compte-bouchaud1998}.
Many mathematical papers followed, a few of which we mention below. The derivation of scaling limits of the models 
is a common theme.

Broadly speaking a trap model is a continuous-time Markov jump process on some regular graph with random transition rates given in terms of 
heavy tailed random variables, the trap environment, which give rise to trapping mechanisms leading to the above mentioned 
effects.
The most studied cases in the mathematics literature involve a jump chain which is independent of the trap environment and spatially
homogeneous, and inverse jump rates given by iid heavy tailed random variables, viewed in this case as trap depths. 
In these cases the trap model 
is thus a time change of the jump chain (which is a discrete time random walk). In this paper we will be concerned with such a trap model 
on $\mathbb Z$, so 
let us discuss the case of $\mathbb Z^d$, $d\geq1$, for a while. (References for the cases of other graphs, like the complete graph 
or the hypercube, may be found in the references mentioned below.)

The simple symmetric case in $d=1$ was studied in~\cite{fin2002}, and a scaling limit was derived, from which aging and localization 
results followed. The higher dimensional symmetric case was resolved in \cite{benarous2006} and \cite{benarouscerny2007}.
In both cases the scaling limit is given by a time change of Brownian motion, with the time change dependent of the Brownian motion 
in $d=1$, but not in $d\geq2$. (The distinction arises as follows: consider the numbers of visits of the jump chain 
to the deepest traps, which in all cases account for virtually all of the time spent by the continuous time process along its history. 
In the first case, this
numbers are macroscopically correlated with the trajectory of the jump chain, and in the limit this manifests itself in the 
representation of the time change in terms of the local time of the scaling limit of the jump chain, as well as in terms of the scaling 
limit of the deep traps. This mechanism is also at play in the model of this paper.
In the second case, those numbers are only weakly correlated with the trajectory of the jump chain, as well as among 
themselves --- the correlations disappear in the scaling limit; this is easy to convince oneself of in the transient case of $d\geq3$, 
but is also the case in the {\em weak} recurrent case of $d=2$. The upshot is that the time change in the limit process is independent
of the scaling limit of the jump  chain.)
Asymptotic aging and localization functions of the trap model are given in terms of the expectations 
of the scaling limit. A variation of this case, is the {\em asymmetric} model, a nearest neighbor model, where the transition rates 
depend on heavy tailed random variables of both origin and destination sites. In this case the jump chain depends on the environment.
Scaling limit and aging results were obtained 
in \cite{benarouscerny2005}, \cite{cerny}, \cite{barlowcerny2011} and \cite{mourrat2011}. The scaling limit is similarly given by 
the time change of a Brownian motion.

Another variation is in the direction of allowing a generic jump chain/ran\-dom walk. This includes the case studied in the present paper. 
Scaling limit and aging results were derived in \cite{fontesandmathieu2013} for the generic case under the validity of a law of large 
numbers for the range of the jump chain, and the slow variation of the tail of the distribution of its return probability. 
These assumptions include all the random walks in $d\geq2$. The process considered in the latter paper is the {\em trap process}, 
namely the depth of the currently visited trap. The scaling limit (which might not exist in the spatial representation of the process) 
is given in terms of a subordinator seen at the inverse of another, correlated subordinator.

In this paper we consider the model on $\mathbb Z$, and assume that the jump chain is a mean zero, $\beta$-stable random walk, with 
$\beta\in(1,2]$, but otherwise generic. This is outside the assumptions of \cite{fontesandmathieu2013}. The model of~\cite{fin2002}, 
where the jump chain is the simple symmetric random walk, is a particular case (of $\beta=2$). One of our motivations is to close a gap 
left by the above papers. (Let us point out that the case where $\beta\in(0,1]$ is included in \cite{fontesandmathieu2013}.) 
We derive the scaling limit of the (spatial version of the) process, given in terms of a time changed $\beta$-stable process, 
and then obtain aging results for the trap model in terms of the scaling limit. 

We call the limit process (given the proper version of the limit heavy tailed random variables) a {\em quasistable process}, 
following the terminology of 
quasidiffusions for the $\beta=2$ case adopted in the literature (see \cite{fin2002} and references therein); see 
also~\cite{kurenok}. Analytical properties 
of quasidiffusions, like the existence and continuity of transition density functions, were crucial in the derivation of 
(non integrated) aging results for the simple symmetric case of \cite{fin2002}. 
The same results can be readily extended to our more general framework for $\beta=2$, but not for $\beta\in(1,2)$, 
where the analogous analytical properties for the corresponding quasistable processes seem to be missing in the literature. 
This point is another of our motivations: to call attention for the lack
of analytical results for a class of processes, namely the quasistable processes, which naturally extends a better known subclass, 
namely the quasidiffusions. 
Without those results, we may nevertheless obtain integrated aging results, if not ordinary aging results. 
(See the following discussion on aging, and Remark~\ref{nonint} below.)

Let us now briefly discuss {\em aging}. Let $\mathcal X_t$ be a generic stochastic process, which might be the trap model 
described above. Consider $Q(s,t)$ a two-time correlation function of $\mathcal X_t$. We call it an {\em aging function}. 
We say that normal aging occurs if there exists a non trivial function $\mathcal Q:\mathbb R^+\rightarrow\mathbb R$ that is the 
limit of $Q(s,t)$ as $t$ and $s$ go to infinity proportionally, that is,
\begin{eqnarray}\label{resnint}
\lim_{\substack{
   t,s\rightarrow\infty \\
   t/s\rightarrow\theta}} 
   Q(s,t)=\mathcal Q(\theta),
\end{eqnarray}   
with $\theta>0$. This is the ordinary, non integrated case, as opposed to the integrated case, where we introduce a random 
time $\mathbb T$ (independent of $\mathcal X_t$), and consider the aging function given by $E[Q(\lambda\mathbb T,\mu\mathbb T)]$, 
with $\mu,\lambda>0$, and the expectation taken with respect to $\mathbb T$. 
We then say that integrated normal aging occurs if there exists a non trivial function 
$\bar{\mathcal Q}: \mathbb R^+\rightarrow\mathbb R$ such that
\begin{eqnarray}\label{resint}
\lim_{\substack{
   \mu,\lambda\rightarrow\infty \\
   \mu/\lambda\rightarrow\theta}} 
   E[Q(\lambda\mathbb T,\mu\mathbb T)]=\bar{\mathcal Q}(\theta),
\end{eqnarray}   
for $\theta>0$.

Typically, $\mathcal Q$ and $\bar{\mathcal Q}$ are decreasing and onto $[0,1]$. In these cases, (aging) results such 
as~(\ref{resnint}) and~(\ref{resint}) may be interpreted as follows, and this
explains the terminology: after observing the process at a large time $t$, the time it takes to 
get a subsequent (reasonably) decorrelated observation is of the order of $t$, indicating that an 
ever increasing slowing down takes place.

The paper is organized as follows. In Section \ref{suprel}, we describe our trap model in detail, and its rescaling,
and state our scaling limit result (Theorem~\ref{mainresult}), proved in Section~\ref{proof}.
Section \ref{intenv} 
is devoted to obtaining an integrated 
aging result 
(Theorem~\ref{bivagingtheorem}) 
for an aging function to be introduced therein. An appendix collects some results on the scaling limit process.\\

A longer version of this paper, including additional results such as convergence of the the trap process,
and a study of localization, 
as well as the simpler case $\alpha\geq1$ (where the scaling limit is an ordinary $\beta$-stable process), 
can be found at {\it http://arxiv.org/pdf/1302.4758.pdf}.

\section{Model and first result}\label{suprel}

Let $\varepsilon=\{\varepsilon_j,\, j\in\mathbb{N}\}$ be a sequence of iid discrete random variables with distribution
function $F$ in the basin of attraction of a stable law with index $\beta\in(1,2]$, such that 
$E(\varepsilon_1)=0$, and $E(e^{it\varepsilon_1})=1$ if and only if $t$ is multiple of $2\pi$.
The latter assumption is well known to imply that the corresponding random walk is aperiodic.

Let now $X=\{X_i,\, i\in\mathbb N\cup\{0\}\}$ be such that $X_0=0$ and for $n\geq1$
\begin{eqnarray}\label{stablerw}
X_n=\sum_{j=1}^n\varepsilon_j.  
\end{eqnarray}
This sequence is called a $\beta$-stable random walk (see \cite{legallrosen1991}); 
this process is also known as long-range random walk.

The object of our study is a continuous time Markov process $\mathcal{X}=\{\mathcal{X}_t: t\geq0\}$ on 
$\mathbb{Z}$ having $X$ as its jump chain,
and whose jump rates are given by $\{\tau_i^{-1}: i\in\mathbb{Z}\}$, 
where $\tau=\{\tau_i: i\in\mathbb{Z}\}$ is a family of iid (strictly) positive random variables in the basin of attraction of a stable law with index $\alpha\in(0,1)$, independent of $X$. Let us point out that the Markov 
property of $\mathcal{X}$ holds for (almost) every fixed realization of $\tau$. The distribution of $\mathcal{X}$
integrated with respect to $\tau$ is {\em not} Markovian.

Our first result is a scaling limit of $\mathcal{X}$. In order to formulate it, we need to introduce
scaling factors, rescaled processes, and limit processes. Let us start with the scaling factors
and the rescaled process.

\paragraph{Scaling factors and rescaled process.}

We recall the well known fact that the assumption on the jump variables $\varepsilon$ implies the following. If $\beta\in(1,2)$, then
there exist constants $c^->0$ and $c^+>0$, and a slowly varying function at infinity $h(\cdot)$ such that
\begin{eqnarray*}
 P(\varepsilon_1<-x)\sim x^{-\beta}(c^-+o(1))h(x)
\end{eqnarray*}
and
\begin{eqnarray*}
 P(\varepsilon_1>x)\sim x^{-\beta}(c^++o(1))h(x),
\end{eqnarray*}
where as usual $f_1(x)\sim f_2(x)$ means $\lim_{x\rightarrow\infty}f_1(x)/ f_2(x)=1$. 
If $\beta=2$, then $H:(0,\infty)\rightarrow(0,\infty)$, with $H(z)=\int_{-z}^zx^2dF(x)$, 
is a slowly varying function at infinity.

It follows that in each case there exists a positive 
slowly varying function $v(\cdot)$ such that as $n\rightarrow\infty$
\begin{eqnarray*}
h(n^{1/\beta}v(n))v^{-\beta}(n)\longrightarrow1, \,\, \mbox{for}\,\, \beta\in(1,2) 
\end{eqnarray*}
and 
\begin{eqnarray*}
H(n^{1/2}v(n))v^{-2}(n)\longrightarrow2(c^++c^-), \,\, \mbox{for}\,\, \beta=2.
\end{eqnarray*}

The assumption on the inverse rate variables $\tau$ implies that
\begin{eqnarray*}
 P(\tau_0>x)\sim x^{-\alpha}(1+o(1))s(x), \,x\geq0,
\end{eqnarray*}
where $s(\cdot)$ is a slowly varying function at infinity.  
It follows that there exists a slowly varying function $w(\cdot)$ such that
\begin{eqnarray*}
s(n^{1/\alpha}w(n))w^{-\alpha}(n)\longrightarrow1
\end{eqnarray*}
as $n\rightarrow\infty$.

Let us define the sequences
\begin{eqnarray}\label{dn}
d_n=n^{1/\beta}v(n),\quad r_n=nd_n^{-1},
\end{eqnarray}
and
\begin{eqnarray}\label{an}
a_n=r_nb_n,
\end{eqnarray}
where
\begin{eqnarray*}\label{bn}
b_n=d_n^{1/\alpha}w(n).
\end{eqnarray*}

We are now ready to define the rescaled process. Let, for $n\geq1$, let
\begin{equation*}
 \mathcal{X}^{(n)}:=\{\mathcal{X}^{(n)}_t= d_n^{-1}\mathcal{X}_{a_nt}\,t\geq0\}.
\end{equation*}

\begin{remark}
 
In the more explicit case where the slowly varying functions entering the distributions of $\varepsilon$ and $\tau$
are asymptotic to constants (say both equal to 1), we get that $a_n=n^{1-\frac1\beta+\frac1{\alpha\beta}}$, 
and $d_n=n^{\frac1\beta}$. By taking $m=a_n$ as scaling factor, we find that $b_n=m^{-\frac1{\beta+\frac1\alpha-1}}$,
and we see a slowing down term of $\frac1\alpha-1$ appearing due to the traps, as compared with the homogeneous
case with no traps, where we would have $b_n=m^{-\frac1{\beta}}$. Except for slowly varying corrections, we have the
same slowing down term in the general case.

\end{remark}

\paragraph{Limit process.}

An ingredient of the limit process is the stable process $Z=(Z_t)_{t\geq0}$ with characteristic function given by 
\begin{eqnarray*}\label{charactstable}
E(e^{isZ_t})=\exp\{-ct|s|^\beta[1+iq\,\mbox{sgn}(s)]\},
\end{eqnarray*}
where $c=-\Gamma(2-\beta)\frac{\cos(\pi\beta/2)}{\beta-1}$ and $q=\frac{c^--c^+}{c^++c^-}\tan(\pi\beta/2)$.

Another ingredient is a bilateral $\alpha$-stable process $V=\{V_x:\, x\in\mathbb R,\, V_0=0\}$ independent of $Z$. 

Let $\phi(t,x)$ be the local time of $Z_t$, that is, let $\phi:\mathbb{R}^+\times\mathbb{R}\to\mathbb{R}^+$ 
be a random function which is jointly continuous with probability one and satisfies
\begin{eqnarray*}
\mathcal L(s: Z_s\in A, \, 0\leq s\leq t)=\int_A \phi(t,x)dx,
\end{eqnarray*}
for any Borel set $A$, where $\mathcal L$ denotes the Lebesgue measure; see \cite{boylan1964}.
Note that it follows from the fact that $Z$ almost surely does not explode at any finite time that
$\phi(t,\cdot)$ is compactly supported for every $t$.
Now define
\begin{equation}\label{eq:S}
 S_t=\int_{-\infty}^\infty \phi(t,x)dV(x),\, t\geq0.
\end{equation} 
The compactnesss of the support of $\phi(t,\cdot)$, as noted right above, and the local finiteness of $V$ 
(as a measure) make~(\ref{eq:S}) well-defined.
From other elementary properties of $\phi$ and $V$, namely $\sup_x\phi(t,x)$ is strictly increasing in $t$ 
and the support of $V$ is the whole line, we get that $S$ is strictly increasing and continuous. 
So it has an ordinary inverse $S^{-1}$.

We are now ready to state our scaling limit result.

\begin{theorem}\label{mainresult}
Let $\{a_n: n\in\mathbb{N}\}$ and $\{d_n: n\in\mathbb{N}\}$ be the sequences defined in (\ref{an}) and (\ref{dn}), respectively. 
We have that 
\begin{equation}\label{eq:main}
\mathcal{X}^{(n)}\Longrightarrow (Z_{S^{-1}_t})_{t\geq0}
\end{equation} 
as $n\rightarrow\infty$, where $\Longrightarrow$ means convergence in distribution on $D([0,\infty),\mathbb R)$ endowed with 
the $J_1$-Skorohod topology.
\end{theorem}

\begin{remark}
Given an arbitrary fixed distribution function $F$, we may replace $V$ by $F$ in~(\ref{eq:S}), and then get a process $Z_{S^{-1}_t}$
(in this generality, $S$ may have flat intervals, so $S^{-1}$ may have to be taken as a generalized inverse).
In this generality, we call $Z_{S^{-1}_t}$ a {\em quasistable process} in an analogy with the term {\em quasidiffusion}, 
used for the case where $Z$ is a Brownian motion. Without this terminology, quasistable processes were introduced in \cite{stone1963}. They are 
strong Markov processes. Some additional properties are stated and proven in an appendix. Unless otherwise mentioned, we will
stick to the $F=V$ case throughout. In the latter context, it is worth emphasizing that the Markov property holds for 
every fixed realization of $V$, and it does not hold for the process integrated with respect to distribution of $V$.

\end{remark}

\begin{remark}
In each side of ``$\Longrightarrow$'' in~(\ref{eq:main}), we have a processes in a {\em random environment}. 
As pointed out in the definition of $\mathcal{X}$ --- see paragraph below~(\ref{stablerw}) ---, 
$\mathcal{X}^{(n)}$ is a Markov process given the {\em environment} $\tau$. 
And in the above remark we have just seen that
$(Z_{S^{-1}_t})_{t\geq0}$ is a Markov process given the {\em environment} $V$. 
The distributions on the right and hand side 
of ``$\Longrightarrow$'' in~(\ref{eq:main}) are those of those processes {\em integrated} with respect to
their respective environments.
\end{remark}

\section{Proof of Theorem~\ref{mainresult}}\label{proof}

\subsection{Preliminaries}\label{prel}

Let $L(n,x)=\sum_{i=0}^n1\{X_i=x\}$ be the local time (occupation time) of the random walk $X$, that is, the number 
of times that $X$ visits the point $x$ up to time $n\in\mathbb N\cup\{0\}$, 
and the rescaled local time and rescaled jump chain
\begin{eqnarray}\label{reslt}
\phi_n(t,x)=r_n^{-1}L([nt],[xd_n]),\,\,\, Z^{(n)}_t=d_n^{-1}X_{[nt]}, 
\end{eqnarray}
for  $t\in[0,\infty)$ and $x\in\mathbb{R}$. It is well-known that the process $(Z^{(n)}_t)_{t\geq0}$ weakly converges on 
$D([0,\infty),\mathbb{R})$ endowed with the $J_1$ topology to $Z$.

\paragraph{Clock process.}

A key element of our analysis is the {\em clock process} associated to $\mathcal{X}$, defined by $C=(C_t)_{t\geq0}$, 
where
\begin{eqnarray}\label{cp}
C_t=\sum_{i=0}^{[t]}\tau_{X_i}T_i, \, t\geq0,
\end{eqnarray}
where $\{T_i:\, i\in\mathbb N\cup\{0\}\}$ is a sequence of iid exponential variables with mean 1 independent of 
$X$ and $\tau$. 

Notice that ${\mathcal X}$ may be represented as $(X_{C^{-1}_t})_{t\geq0}$, where $C^{-1}$ is the generalized (right continuous)
inverse of $C$.

We have that the clock process (\ref{cp}) is equal in distribution to the process $\bar C=(\bar C_t)_{t\geq0}$,
where
\begin{eqnarray*}\label{ecp}
\bar{C}_t=\sum_{i\in\mathbb{Z}}\tau_i\sum_{j=1}^{L([t],i)}E_{ij},
\end{eqnarray*}
and $\mathcal{E}=\{E_{ij}: i\in\mathbb{Z}, j\in\mathbb{N}\}$ is a family of iid exponential random variables 
with mean 1 and independent of all random variables defined previously; we here define $\sum_{j=1}^0E_{ij}=0$.

We thus have that $(X_{\bar C^{-1}_t})_{t\geq0}$ is a version of $\mathcal{X}$. Furthermore, defining
\begin{eqnarray}\label{scp}
\bar{C}^{(n)}_t\equiv a_n^{-1}\bar C_{nt},\,t\geq0,
\end{eqnarray}
$n\geq0$, we may check that 
\begin{eqnarray*}\label{storep}
\big\{Z^{(n)}_{\bar{C}^{(n)^{-1}}_t}: t\geq0\big\}\stackrel{d}=\mathcal{X}^{(n)},
\end{eqnarray*}
where ``$\stackrel{d}=$'' means equality in distribution.

\paragraph{A convenient version of $\tau$.}

The proof of Theorem~\ref{mainresult} will involve obtaining the scaling limit of the clock process.
Following a strategy used numerously before (for an early reference, see~\cite{fin2002}, Section 3), 
we will to resort to a version of the rescaled trap environment which converges strongly
(rather than only weakly, which is the case for the original trap environment). 
With the new version of $\tau$, we define a new version of the clock process, the rescaling of which we 
then will later on show converges in distribution for almost every realization of the trap environment.
This convergence is a key ingredient of our proof of the scaling limit of $\mathcal X_t$. 

We now present the new version of the environment.
Let $V=\{V_x: x\in\mathbb{R}\}$ be a bilateral $\alpha$-stable process independent of $\mathcal E$ and $X$. 
Consider a function $G:[0,\infty)\rightarrow[0,\infty)$ satisfying 
$ P(V_1>G(y))= P(\tau_0>y)$, $y>0$,
and for $n\geq0$ let $g_n:\mathbb{R}^+\rightarrow[0,\infty)$ be such that 
$g_n(y)=b_n^{-1}G^{-1}(d_n^{1/\alpha}y)$. For $n\geq0$, let
\begin{eqnarray*}\label{taun}
\tau_x^{(n)}\equiv b_ng_n(V_{x+d_n^{-1}}-V_x),\quad x\in d_n^{-1}\mathbb{Z}.
\end{eqnarray*}

One readily checks that $\tau^{(n)}=\{\tau_x^{(n)}: x\in d_n^{-1}\mathbb{Z}\}$ has the same distribution as
$\{\tau_i: i\in \mathbb{Z}\}$ for every $n$. It follows that the process (\ref{scp}) follows the same  
law as that of the following process:
\begin{eqnarray}\label{fv}
\widetilde{C}^{(n)}_t\equiv\sum_{x\in d_n^{-1}\mathbb{Z}}g_n(V_{x+d_n^{-1}}-V_x)\phi_n(t,x)\bar{T}_{xd_n}(nt), \quad t\geq0,
\end{eqnarray}
where for every $t$ and $i$
\begin{equation*}
 \bar{T}_{i}(t)=
\begin{cases}
 \dfrac{\sum_{j=1}^{L([t],i)}E_{ij}}{L([t],i)},&\mbox{ if } L([t],i)>0;\\
 \quad\quad\quad\quad 1,&\mbox{ otherwise}, 
\end{cases}
\end{equation*} 
and $\phi_n(t,x)$ is the rescaled local time defined in (\ref{reslt}).

We thus have that for every $n$
\begin{eqnarray*}\label{storep1}
\big\{Z^{(n)}_{\widetilde{C}^{(n)^{-1}}_t}: t\geq0\big\}\stackrel{d}=\mathcal{X}^{(n)}.
\end{eqnarray*}

We will use the left hand side version of $\mathcal{X}^{(n)}$ in the proof of Theorem~\ref{mainresult}
next.

\subsection{Proof of Theorem~\ref{mainresult}}\label{lim}

\paragraph{Path space and topology.}

Denote by $D([0,T],\mathbb R)$ and $D([0,T],\mathbb R^+)$ the space of the real and non negative c\`adl\`ag functions 
on $[0,T]$, where $T>0$, and
let $J_{1T}$ and $u_T^+$ be the $J_1$-Skorohod and uniform metrics in $D([0,T],\mathbb R)$ and $D([0,T],\mathbb R^+)$, 
respectively. Further let $d=\sum_{n=1}^\infty2^{-n}\min(J_{1n},1)$ and $u^+=\sum_{n=1}^\infty2^{-n}\min(u_n^+,1)$. 
We denote the uniform topology by $U$.

\paragraph{An auxiliary result.}

In order to prove the Theorem \ref{mainresult}, we first obtain the joint scaling limit of the clock process and the jump chain 
$X$, provided in the next result. This strategy was used, for example, in \cite{benarouscerny2007} to find the scaling limit
for the trap model on $\mathbb Z^d$ (for $d\geq2$) 
with nearest neighbors and inverse rates as we consider here. The major work lies in showing the joint convergence of the 
finite-dimensional distributions of the clock process and the jump chain. 
In \cite{benarouscerny2007} analytical arguments for the joint characteristic function was used. 
We here consider a different approach based on a probabilistic argument. One important difference between our clock process and that 
of \cite{benarouscerny2007} is that in our case the scaling limit depends on the scaling limit of the rescaled random walk $Z_t^{(n)}$ 
(in \cite{benarouscerny2007} the scaling limits of the clock process and the jump chain are independent). 

\begin{theorem}\label{bivweaklimit}
For almost every realization of $V$, 
$(Z^{(n)},\widetilde C^{(n)})=\{(Z^{(n)}_t,\widetilde{C}^{(n)}_t): t\geq0\}$ converges to 
$(Z,S)=\{(Z_t,S_t): t\geq0\}$ in distribution as $n\rightarrow\infty$ on 
$D([0,\infty),\mathbb R)\times D([0,\infty),\mathbb R^+)$ endowed with 
the $J_1\times U$ product topology.
\end{theorem}

\begin{remark}\label{contsl}
The limit process of our rescaled clock process is almost surely continuous while the one considered in 
\cite{benarouscerny2007} is not (there, the limit clock process is an $\alpha$-stable subordinator). 
Another difference between both cases is with respect to the topology where the convergence takes place. 
It holds with the uniform topology here, while in \cite{benarouscerny2007} one needs to use the $M_1$-Skorohod 
topology. Nevertheless, our proof will indeed verify criteria of convergence in the $J_1\times M_1$ topology for the 
bivariate process given in Theorem~\ref{bivweaklimit}. Since $S$ is almost surely continuous, the convergence in 
$J_1\times U$ topology follows --- see e.g.~\cite{whitt2002}, Subsection 3.3.
\end{remark}

\begin{remark}
From now on we will denote the conditional distribution and expectation given $V$ by $\mathbb P(\cdot)\equiv P(\cdot|V)$ 
and $\mathbb E(\cdot)\equiv E(\cdot|V)$, respectively.
\end{remark}

\noindent {\bf Proof of Theorem \ref{bivweaklimit}.} 
By a standard argument, it is enough to show convergence on $D([0,T],\mathbb R)$ and $D([0,T],\mathbb R^+)$ 
equipped with the metrics $J_{1T}$ and $u_T^+$, respectively, with $T>0$ arbitrary. 
We will next prove the convergence of the finite-dimensional distributions of $(Z^{(n)},\widetilde{C}^{(n)})$ 
and then establish tightness in those path spaces. This then implies the result.\\

\noindent\emph{Convergence of the finite-dimensional distributions of $(Z^{(n)},\widetilde C^{(n)})$}\\

We resort to a result by
Borodin (1985), which, under conditions ensured by our assumptions on $\varepsilon$,
guarantees the existence on some probability space of processes
${Z'}^{(n)}=\{{Z_t'}^{(n)}: t\in[0,T]\}_{n\geq1}$ and $Z'=\{Z'_t: t\in[0,T]\}$ such that
\begin{itemize}
 \item[(a)] their finite-dimensional distributions coincide with those of $Z^{(n)}$ and $Z$, 
respectively;
 \item[(b)] ${Z'}^{(n)}$ converges almost surely to $Z'$ on $D([0,T],\mathbb R)$ endowed with the $J_1$-Skorohod topology;
 \item[(c)] the local times $\phi_n'(\cdot,\cdot)$ and $\phi'(\cdot,\cdot)$ of ${Z'}^{(n)}$ and $Z'$, 
respectively, are such that for any $T>0$ and $\xi>0$:
\begin{eqnarray}\label{convlocaltime}
\lim_{n\rightarrow\infty}P\big(\sup_{(t,x)\in[0,T]\times\mathbb{R}}|\phi'_n(t,x)-\phi'(t,x)|>\xi\big)=0
\end{eqnarray}
\end{itemize}
(see Theorem 1.1 in~\cite{borodin1985}).

Therefore, it is enough to 
show the convergence of the finite-dimensional distributions of 
$({Z'}^{(n)},{\widetilde{C}}'^{(n)})$ to 
those of $(Z,S)$,
where ${\widetilde{C}}'^{(n)}$ is defined analogously as (\ref{fv}), 
replacing the quantities that depend on $Z^{(n)}$ by the corresponding ones depending on ${Z'}^{(n)}$. 
We will likewise below use the notation $\mathcal B'$ when replacing $Z^{(n)}$ by ${Z'}^{(n)}$ in a quantity
$\mathcal B$ depending on the former process.

We now define the sets of {\em deep traps}. For arbitrary $\delta>0$, let
\begin{eqnarray*}{}
\mathcal{T}_\delta&=&\{x\in\mathbb{R}: V_x-V_{x-}>\delta\}=\{\ldots<x_{-1}<x_0<x_1<\ldots\},\\
\mathcal{T}_\delta^{(n)}&=&\{x_j^{(n)}: j\in\mathbb{Z}\},\,n\geq1,
\end{eqnarray*}
where
$x_{j}^{(n)}=d_n^{-1}[d_nx_j]$, $j\in\mathbb{Z}$.
We first show that 
\begin{eqnarray}\label{trapc}
\lim_{\delta\downarrow0}\limsup_{n\rightarrow\infty}\sum_{x\in d_n^{-1}\mathbb{Z}\cap(\mathcal{T}_\delta^{(n)})^c}
g_n(V_{x+d_n^{-1}}-V_x)\phi_n'(t,x)\bar{T}_{xd_n}'(nt)=0,
\end{eqnarray}
in probability, for all $t\in[0,T]$. This says that the terms of the rescaled clock process that are out of the deep trap set 
have a negligible contribution to the limit process as $n\rightarrow\infty$ and $\delta\downarrow0$.

\begin{remark}
Let $\mathcal A^{(n,\delta)}_t$ be a random variable depending on $\delta$, $n$ and $t$. 
Below, when we say that $\lim_{\delta\downarrow0}\limsup_{n\rightarrow\infty} \mathcal A^{(n,\delta)}_t=0$ in probability, 
that
means that for all $\epsilon>0$ fixed we have that 
$\lim_{\delta\downarrow0}\limsup_{n\rightarrow\infty}P(|\mathcal A^{(n,\delta)}_t|>\epsilon)=0$.
\end{remark}

Borodin's result mentioned above then says that for any $\epsilon>0$, there exists $A_{\epsilon}>0$ such that 
\begin{eqnarray}\label{borodintheorem}
\sup_n P\bigg(\sup_{t\in[0,T]}|{Z_t'}^{(n)}|>A_\epsilon\bigg)<\epsilon\quad\mbox{and}
\quad P\bigg(\sup_{t\in[0,T]}|Z_t'|>A_\epsilon\bigg)<\epsilon.
\end{eqnarray}

Let $I_\epsilon=(-A_\epsilon, A_\epsilon)$. Then, for all $\zeta>0$, using the above result and the Markov inequality, 
we get that
\begin{eqnarray}\label{auxsum}
&&\mathbb P\bigg(\sum_{x\in d_n^{-1}\mathbb{Z}\cap(\mathcal{T}_\delta^{(n)})^c}
g_n(V_{x+d_n^{-1}}-V_x)\phi_n'(t,x)\bar{T}_{xd_n}'(nt)>\zeta\bigg)\nonumber\\
&\leq&\mathbb P\bigg(\sum_{x\in d_n^{-1}\mathbb{Z}\cap(\mathcal{T}_\delta^{(n)})^c\cap I_\epsilon} 
g_n(V_{x+d_n^{-1}}-V_x)\phi_n'(t,x)\bar{T}_{xd_n}'(nt)>\zeta\bigg)+\epsilon\nonumber\\
&\leq&\zeta^{-1}\sum_{x\in d_n^{-1}\mathbb{Z}\cap(\mathcal{T}_\delta^{(n)})^c\cap I_\epsilon}
g_n(V_{x+d_n^{-1}}-V_x)E(\phi_n'(t,x))+\epsilon
\end{eqnarray}  
for almost every realization of $V$. For all $M>0$ integer we have that
\begin{eqnarray}\label{moments}
E(\phi_n^M(t,x))\leq E(\phi_n^M(t,0)),\,
\lim_{n\rightarrow\infty}\!\!E(\phi_n^M(t,0))=\dfrac{t^{M(1-1/\beta)}z\Gamma(1-1/\beta)}{\Gamma(2-1/\beta)},
\end{eqnarray}
where $z$ is the value of the density of $Z_1$ at zero (see \cite{borodin1985}, page 328). 
From (\ref{moments}) and using the equality of the finite-dimensional distributions of $\phi_n(t,x)$ and $\phi_n'(t,x)$, 
which holds for all $x\in\mathbb R$, we obtain that the sum in (\ref{auxsum}) is bounded above by constant times
$$\sum_{x\in d_n^{-1}\mathbb{Z}\cap(\mathcal{T}_\delta^{(n)})^c\cap I_\epsilon}g_n(V_{x+d_n^{-1}}-V_x).$$

Now, arguing as in \cite{fin2002}, paragraphs of (3.25) to (3.28), it follows that 
$\lim_{\delta\downarrow0}\limsup_{n\rightarrow\infty}$ of the above term vanishes in probability. 
This and the arbitrariness of $\epsilon$ yield (\ref{trapc}).  

With the above result, we now define the clock process restricted to the deep traps:
\begin{eqnarray}\label{relprocess}
{\widetilde{C}_t}'^{(n,\delta)}=\sum_{i\geq1}g_n(V_{x_i^{(n)}+d_n^{-1}}-V_{x_i^{(n)}})\phi_n'(t,x_i^{(n)})\bar{T}_{x_i^{(n)}d_n}'(nt).
\end{eqnarray}

Let $\epsilon>0$ and take $A_\epsilon$ satisfying (\ref{borodintheorem}). Define the set 
$$\Omega_{\epsilon,n}=\{\sup_{t\in[0,T]}|Z^{(i)}_t|\leq A_\epsilon, 1\leq i\leq n\}\cap\{\sup_{t\in[0,T]}|Z_t|\leq A_\epsilon\}.$$

On $\Omega_{\epsilon,n}$, the process given in (\ref{relprocess}) equals
\begin{eqnarray*}\label{relprocess2}
\widetilde{C}'^{(n,\delta,\epsilon)}_t=\sum_{i=-N_{\delta,\epsilon}}^{N_{\delta,\epsilon}}
g_n(V_{x_i^{(n)}+d_n^{-1}}-V_{x_i^{(n)}})\phi_n'(t,x_i^{(n)})\bar{T}_{x_i^{(n)}d_n}'(nt),
\end{eqnarray*}
where $N_{\delta,\epsilon}=\max\{|j|\in\mathbb{N}: x_j\in I_{\epsilon}\}<\infty$ for almost every $V$. 

Result (\ref{convlocaltime}) and the a.s.~continuity of $\phi'(\cdot,\cdot)$ 
imply that $\phi_n'(t,x_i^{(n)})$ converges in probability to $\phi'(t,x_i)$ 
uniformly in $i\in[-N_{\delta,\epsilon},N_{\delta,\epsilon}]$ and $t\in[0,T]$ as $n\rightarrow\infty$. 
Furthermore, since $L'([nt],[x_i^{(n)}d_n])\stackrel{p}{\longrightarrow}\infty$ when $n\rightarrow\infty$ 
(this follows from the convergence in probability of $\phi_n'(t,x_i^{(n)})$), it follows from the Law of the 
Large Numbers that $\bar{T}'_{x_i^{(n)}d_n}(nt)\stackrel{p}{\longrightarrow}1$ as $n\rightarrow\infty$ for all 
$i$ and $t$.
Further, from Proposition 3.1 of \cite{fin2002}, it follows that 
\begin{eqnarray}\label{prop3.1fin2002}
g_n(V_{x_i^{(n)}+d_n^{-1}}-V_{x_i^{(n)}})\longrightarrow V_{x_i}-V_{x_i-},
\end{eqnarray}
as $n\rightarrow\infty$, for almost every $V$ and for all $i$.

The convergence in probability of $\phi_n'(t,x_i^{(n)})$ and $\bar{T}'_{x_i^{(n)}d_n}(nt)$ 
(uniformly in $i\in[-N_{\delta,\epsilon},N_{\delta,\epsilon}]$) discussed above, result~(\ref{prop3.1fin2002}),
and the Mapping Theorem imply that for almost every $V$  
\begin{eqnarray}\label{cprob1}
\big({Z_t'}^{(n)},{{{\widetilde{C}}_t}}^{\prime(n,\delta,\epsilon)}\big)\stackrel{p}{\longrightarrow}\big(Z'_t,{S_t'}^{(\delta,\epsilon)}\big),
\end{eqnarray}
as $n\rightarrow\infty$ for every $t$, where
${S_t'}^{(\delta,\epsilon)}=\sum_{x\in\mathcal{T}_\delta\cap I_\epsilon}(V_x-V_{x-})\phi'(t,x)$. 
Moreover, we have that 
\begin{eqnarray}\label{cprob2}
\lim_{\delta\downarrow0}\lim_{\epsilon\downarrow0}{S_t'}^{(\delta,\epsilon)}=\int_{-\infty}^\infty\phi'(t,x)dV_x,
\end{eqnarray}
almost surely. The convergence of the finite-dimensional distributions of 
$(Z^{(n)},$ $\widetilde{C}^{(n)})$ to those of $(Z,S)$ now follows from (\ref{cprob1}) and (\ref{cprob2}).\\

\noindent\emph{Tightness}\\

We now show that the sequence $(Z^{(n)},\widetilde{C}^{(n)})$ is tight on 
$D([0,T],\mathbb R)\times D([0,T],\mathbb R^+)$ endowed with the $J_1\times M_1$ product topology. 
It is enough to establish tightness
of each coordinate. The first coordinate converges in distribution, so it is tight. We are thus left with showing that 
$(\widetilde{C}^{(n)})$ is tight on $D([0,T],\mathbb R^+)$ endowed with the $M_1$ topology (and consequently with the 
$U$ topology --- see Remark~\ref{contsl}). 

Using the fact that $\widetilde{C}^{(n)}$ is non decreasing, one readily checks that the oscillation function 
$w_s$ used in the condition (ii) of the Theorem 12.12.3 of \cite{whitt2002} equals 0. With this and using the fact 
that the process is non negative, we have that the tightness criteria of the latter theorem are given by
\begin{eqnarray*}
&&\mbox{(i)} \lim_{c\rightarrow\infty}\limsup_{n\rightarrow\infty}\mathbb P\big(\widetilde{C}^{(n)}_T>c\big)=0;\\
&&\mbox{(ii) For each $\xi>0$},\, \lim_{\varepsilon\downarrow0}\limsup_{n\rightarrow\infty}
\mathbb P\big(\max\{\bar{v}(C^*_n,0,\varepsilon),\bar{v}(C^*_n,T,\varepsilon)\}>\xi\big)=0,
\end{eqnarray*}
where, for $t\in[0,T]$, we have that $$\bar{v}(C^*_n,t,\varepsilon)
=\sup_{\max\{0,t-\varepsilon\}\leq t_1\leq t_2\leq \min\{t+\varepsilon,T\}}
\{\widetilde{C}^{(n)}_{t_2}-\widetilde{C}^{(n)}_{t_1}\}.$$

Taking $\varepsilon<T$, we have that, for $t=0$ and $t=T$, the quantity above reduces to 
\begin{eqnarray}\label{v1}
\bar{v}(C^*_n,0,\varepsilon)=\widetilde{C}^{(n)}_\varepsilon
\end{eqnarray}
and
\begin{eqnarray}\label{v2}
\bar{v}(C^*_n,T,\varepsilon)=\widetilde{C}^{(n)}_T-\widetilde{C}^{(n)}_{T-\varepsilon}.
\end{eqnarray}
respectively. 

We now show that the rescaled clock process satisfies conditions (i) and (ii). 
From now on, we use the restriction $\varepsilon<T$. Condition (i) follows from the convergence of 
the finite-dimensional distributions, that is
\begin{eqnarray*}
\lim_{c\rightarrow\infty}\limsup_{n\rightarrow\infty}\mathbb P\big(\widetilde{C}^{(n)}_T>c\big)=\lim_{c\rightarrow\infty}\mathbb P\big(S_T>c\big)=0.
\end{eqnarray*}

Using (\ref{v1}) and the convergence of the finite-dimensional distributions, we also get that for each $\xi>0$ fixed
$$\lim_{\varepsilon\downarrow0}\limsup_{n\rightarrow\infty}\mathbb P\big(\bar{v}(C^*_n,0,\varepsilon)>\xi\big)=\lim_{\varepsilon\downarrow0}
\mathbb P(S_\varepsilon>\xi).$$

Since $S_\varepsilon$ converges to 0 in probability as $\varepsilon\downarrow0$, we find that
\begin{eqnarray}\label{limv1}
\lim_{\varepsilon\downarrow0}\limsup_{n\rightarrow\infty}\mathbb P\big(\bar{v}(C^*_n,0,\varepsilon)>\xi\big)=0.
\end{eqnarray}

Similarly as in the previous case, using (\ref{v2}), for each $\xi>0$ fixed, we have that
$$\lim_{\varepsilon\downarrow0}\limsup_{n\rightarrow\infty}
\mathbb P\big(\bar{v}(C^*_n,T,\varepsilon)>\xi\big)=\lim_{\varepsilon\downarrow0}\mathbb P(S_T-S_{T-\varepsilon}>\xi).$$

Now using the fact that $S$ is almost surely continuous, we get that
$S_T-S_{T-\varepsilon}\stackrel{a.s.}{\longrightarrow}0$ as $\varepsilon\downarrow0$.
With this, we obtain that
\begin{eqnarray}\label{limv2}
\lim_{\varepsilon\downarrow0}\limsup_{n\rightarrow\infty}\mathbb P\big(\bar{v}(C^*_n,T,\varepsilon)>\xi\big)=0.
\end{eqnarray}

Results (\ref{limv1}) and (\ref{limv2}) imply that the condition (ii) is satisfied, and hence 
$\widetilde{C}^{(n)}$ is tight.  $\square$\\

\noindent {\bf Proof of Theorem \ref{mainresult}.} 
Let us start by defining some subsets of $D([0,\infty),\mathbb R)$ (definitions which hold 
analogously in the case of $D([0,\infty),\mathbb R^+)$). 

Let $D_\uparrow([0,\infty),\mathbb R)$ be the space of the non negative functions in $D([0,\infty),\mathbb R)$ 
that are non decreasing. We here denote the space of continuous functions which are strictly increasing by 
$C_{\uparrow\uparrow}([0,\infty),\mathbb R)$. Denote by $D_u([0,\infty),\mathbb R)$ be the space of functions in 
$D([0,\infty),\mathbb R)$ which are unbounded. Hence, we define $D_{u,\uparrow}([0,\infty),\mathbb R)
=D_u([0,\infty),\mathbb R)\cap D_\uparrow([0,\infty),\mathbb R)$.

Using Theorem \ref{bivweaklimit} and the arguments in its proof, one may check that the 
finite-dimensional distributions of 
$(Z^{(n)},\widetilde{C}^{(n)^{-1}})$ weakly converge to those of $(Z,S^{-1})$ as $n\rightarrow\infty$. 

Further, we have that $D_{u,\uparrow}([0,\infty),\mathbb R^+)$ endowed with the uniform topology is separable 
and complete, so by using the weak convergence of $\widetilde{C}^{(n)^{-1}}_t$ to $S_t^{-1}$ 
(which follows from Theorem \ref{bivweaklimit} and the Mapping Theorem) and the converse half of Prohorov's 
Theorem, we get tightness of $(\widetilde{C}^{(n)^{-1}})$. 
Therefore $(Z^{(n)},\widetilde{C}^{(n)^{-1}})\Longrightarrow(Z,S^{-1})$ as $n\rightarrow\infty$ 
on $D([0,\infty),\mathbb R)\times D_{u,\uparrow}([0,\infty),\mathbb R^+)$ equipped with the $J_1\times U$ 
product topology.

Now we use Theorem 13.2.2 of \cite{whitt2002}, which states that the composition map from 
$D([0,\infty),\mathbb R)\times D_{u,\uparrow}([0,\infty),\mathbb R^+)$ to $D([0,\infty),\mathbb R)$ is 
continuous on $D([0,\infty),\mathbb R)\times C_{\uparrow\uparrow}([0,\infty),\mathbb R^+)$ equipped with 
the $J_1$ topology. Notice that the trajectories of 
$(Z,S^{-1})$ are in $D([0,\infty),\mathbb R)\times C_{\uparrow\uparrow}([0,\infty),\mathbb R^+)$
almost surely. 
The Mapping Theorem and the above results imply that $\mathcal{X}^{(n)}\Longrightarrow Z_{S^{-1}}$ on 
$D([0,\infty),\mathbb R)$ endowed with the $J_1$, yielding the result. $\square$\\ 

We close this section with a result about the self-similarity of the scaling limit. This will prove 
useful to arguing our aging result. We say that a process 
$\{W(t): t\geq0\}$ is self-similar of order $H>0$ if $W(at)\stackrel{fdd}{=}a^HW(t)$, for all $a>0$, where 
``$\stackrel{fdd}{=}$'' means equality of all finite-dimensional distributions.

\begin{proposition}\label{selfsimilar} The process $\{Z_{S_t^{-1}}: t\geq0\}$ is self-similar of order 
$\alpha\beta^2/\{1+\alpha(\beta-1)\}$. \end{proposition}

\noindent {\bf Proof.} We first claim that $\{S_t: t\geq0\}$ is self-similar of order $1-1/\beta+1/(\alpha\beta)$. 
Let $a>0$. By using the self-similarity property of the process $Z_t$ (of order $1/\beta$) and the 
following representation for the local time of the process $Z_t$: 
$$\phi(t,x)=\lim_{\epsilon\downarrow0}\frac{1}{2\epsilon}\int_0^t1\{x-\epsilon<Z_s\leq x+\epsilon\}ds,$$
which holds for $t>0$ and $x\in\mathbb{R}$, one readily checks that 
$$\phi(at,x)\stackrel{fdd}{=}a^{1-1/\beta}\phi(t,a^{-1/\beta}x).$$ 

With this and using the self-similarity (of order $1/\alpha$) of the process $V_t$, 
we readily justify our first claim by making a straightforward change of variable in the integral defining 
of $S$ (see~(\ref{eq:S})).
Hence we have that $$S^{-1}_{at}\stackrel{fdd}{=}a^{1/\{1-1/\beta+1/(\alpha\beta)\}}S^{-1}_{t}.$$ 

This and the self-similarity of $Z$ imply that $Z_{S_t^{-1}}$ is self-similar of order 
$\beta/\{1-1/\beta+1/(\alpha\beta)\}=\alpha\beta^2/\{1+\alpha(\beta-1)\}$. $\square$\\

\section{Integrated aging}\label{intenv}

In the final section we move our attention to the study of aging for $\mathcal X$. 
We prove an integrated aging result, as explained at the introduction. 
Consider the following integrated aging function for $\mathcal{X}_t$: 
\begin{eqnarray}\label{intaging}
\bar{\mathcal{R}}(\lambda,\mu)=E[R(\lambda\mathbb{T},\mu\mathbb{T})],
\end{eqnarray}
where $\mu,\lambda\geq0$, $\mathbb{T}$ is an absolutely continuous random variable 
supported on $(0,\infty)$ and independent of all other variables, and
\begin{eqnarray*}
R(s,t)=P(\mathcal{X}_t=\mathcal{X}_{t+s}),
\end{eqnarray*}
for $s,t>0$. 

We now state and prove an integrated aging result for (\ref{intaging}). 

\begin{theorem}\label{bivagingtheorem} Let $\mathcal{R}: \mathbb{R}^+\rightarrow[0,1]$ such that 
$\mathcal{R}(\theta)=P(Z_{S_1^{-1}}=Z_{S_{1+\theta}^{-1}})$. Then  
\begin{eqnarray*}\label{bivaging}
\lim_{\substack{
   \lambda,\mu\rightarrow\infty \\
   \lambda/\mu\rightarrow\theta}} 
   \bar{\mathcal{R}}(\lambda,\mu)=\mathcal{R}(\theta).
\end{eqnarray*}
\end{theorem}

\begin{remark}
In the appendix we give an argument for the non triviality of $\mathcal R(\cdot)$. (See Corollary~\ref{coro}.)
\end{remark}

\begin{remark}\label{fin}
In \cite{fin2002}, non integrated aging results were established for a number of aging functions of $\mathcal X$, 
including the one of Theorem \ref{bivagingtheorem}, for the case where $X$ is the simple symmetric random walk in dimension 1. 
Recall the discussion at the introduction on the limitations on extending this approach to the more general situation of
the present paper.
\end{remark}

The convergence in distribution stated in Theorem~\ref{mainresult} is not sufficient to prove Theorem \ref{bivagingtheorem}. 
Additional arguments will be stated and proven in two lemmas which will allow us to replace the process $\mathcal X_t$ 
by a process living in the deep traps. 
This will lead to the desired result. 

For simplicity of notation, we define 
$$Y_t^{(n)}=Z^{(n)}_{\widetilde{C}^{(n)^{-1}}_t},\,\, Y_t^{(n,\delta)} =Z^{(n)}_{\widetilde{C}^{(n,\delta)^{-1}}_t},
\,\, Y_t^{(\delta)}=Z_{S^{(\delta)^{-1}}_t},\,\, t\geq0,$$
where ${S_t^{(\delta)}}^{-1}$ is the inverse of 
$S^{(\delta)}_t=\sum_{x\in\mathcal{T_{\delta}}}\phi(t,x)(V_x-V_{x-})$. 

We now introduce an auxiliary process, denoted by 
$\bar{Y}_t^{(n,\delta)}$, which lives on $\delta$-traps. To define it, let $W_i^{(n,\delta)}$, $i=0,1,\ldots$, 
be the successive $\delta$-traps visited by $Y_t^{(n)}$, with the restriction that $W_i^{(n,\delta)}\neq W_{i+1}^{(n,\delta)}$, 
and let $\bar S_i^{(n,\delta)}$, $i=0,1,\ldots$, denote the successive hitting times of those traps by  $Y_t^{(n)}$, respectively
(so that $Y_{\bar S_i^{(n,\delta)}}^{(n)}=W_i^{(n,\delta)};\,Y_{\bar S_i^{(n,\delta)}-}^{(n)}\ne W_i^{(n,\delta)}$).
Let us make
\begin{eqnarray*}
\bar{Y}_t^{(n,\delta)}=W_{i}^{(n,\delta)}, \quad \mbox{if} \quad \bar S_i^{(n,\delta)}\leq t< \bar S_{i+1}^{(n,\delta)}, \quad i=0,1,\ldots
\end{eqnarray*}

Notice that $Y_t^{(n)}=\bar{Y}_t^{(n,\delta)}$ whenever $Y^{(n)}$ is visiting a $\delta$-trap, and different otherwise,
and $\bar{Y}^{(n,\delta)}$ of course lives on $\delta$-traps.

\begin{lemma}\label{lemma1aging}
For any $T>0$, we have that 
$$\lim_{\delta\downarrow0}\limsup_{n\rightarrow\infty}\int_0^T \mathbb P(Y_t^{(n)}\neq \bar{Y}_t^{(n,\delta)})\,dt=0,$$
for almost every $V$.
\end{lemma}

\begin{remark}\label{nonint}
A strengthening of this lemma to a non integrated version would lead to a non integrated version of 
Theorem~\ref{bivagingtheorem}. As pointed out at the introduction, and in comparison with the 
approach of~\cite{fin2002} to obtaining non integrated aging results, we lack here analogues of the analytical
results for quasistable processes which exist for quasidiffusions.
\end{remark}

\noindent {\bf Proof of Lemma~\ref{lemma1aging}.} Since $\phi(t,x)\stackrel{a.s.}{\longrightarrow}\infty$ as $t\rightarrow\infty$, 
for any $x\in\mathbb{R}$ --- see for instance \cite{boylan1964} ---, we have that 
$S_t^{(\delta)}\stackrel{a.s.}{\longrightarrow}\infty$ as $t\rightarrow\infty$, for all fixed $\delta>0$.
This and the weak convergence  
$\lim_{\delta\downarrow0}\lim_{n\rightarrow\infty}\widetilde C_t^{(n,\delta)}\stackrel{d}{=}S_t$, 
for each $t$ (which can be easily obtained from the elements of the proof of Theorem~\ref{bivweaklimit})
we get that, given $T,\eta,\delta>0$ there exist $n_0,\Delta>0$ such that 
\begin{eqnarray}\label{S:clock}
\mathbb P(\widetilde C_\Delta^{(n)}\leq T)\leq \mathbb P(\widetilde C_\Delta^{(n,\delta)}\leq T)\leq\eta, \quad \forall n\geq n_0.
\end{eqnarray}
We then fix $T,\eta,\delta>0$, take $n_0,\Delta>0$ such that the above inequalities hold and obtain that
\begin{eqnarray*}
&&\int_0^T \mathbb P(Y_t^{(n)}\neq \bar{Y}_t^{(n,\delta)})dt
=\mathbb E\left(\int_0^T 1\{Y_t^{(n)}\neq \bar{Y}_t^{(n,\delta)}\}dt\right)\\
&\leq&\mathbb E\left(\int_0^T 1\{Y_t^{(n)}\neq \bar{Y}_t^{(n,\delta)},\widetilde C_\Delta^{(n,\delta)}\geq T\}dt\right)+\eta T\\
&\leq&\mathbb E\bigg(\sum_{x\in d_n^{-1}\mathbb{Z}\cap(\mathcal{T}_\delta^{(n)})^c}
g_n(V_{x+d_n^{-1}}-V_x)\phi_n(\Delta,x)\bar{T}_{xd_n}(n\Delta)\bigg)+\eta T,
\end{eqnarray*}
The above inequality, (\ref{trapc}) and the arbitrariness of $\eta$ yield the result. $\square$ \\

\begin{lemma}\label{lemma2aging}
For almost every realization of $V$ we have that 
$$\lim_{\delta\downarrow0}\limsup_{n\rightarrow\infty}J_1((\bar{Y}_t^{(n,\delta)}),(Y_t^{(n,\delta)}))=0,$$
in probability, where $J_1$ is the Skorohod metric.
\end{lemma}
\noindent {\bf Proof.} The processes 
$\bar Y^{(n,\delta)}=\{\bar Y_t^{(n,\delta)}: t\geq0\}$ and $Y^{(n,\delta)}=\{Y_t^{(n,\delta)}: t\geq0\}$ 
successively visit the same traps but with different sojourn times. 
So, it is enough to show that the maximum of the differences between the successive sojourn times of the 
traps visited by both processes within $[0,T]$, respectively,  goes to 0 in probability as $n\rightarrow\infty$ 
and $\delta\downarrow0$. 
Let $S_i^{(n,\delta)}$ be the successive jump times of $Y^{(n,\delta)}$.
We then have that 
$$Y_t^{(n,\delta)}=W_i^{(n,\delta)}, \,\,\mbox{if}\,\,  S_i^{(n,\delta)}\leq t<  S_{i+1}^{(n,\delta)},$$
for $i\in\mathbb{N}\cup\{0\}$.

Let $\mathcal{S}_i^{(n,\delta)}=S_i^{(n,\delta)}-S_{i-1}^{(n,\delta)}$, 
$\bar{\mathcal{S}}_i^{(n,\delta)}=\bar S_i^{(n,\delta)}-\bar S_{i-1}^{(n,\delta)}$, 
$i\geq1$, denote the successive
sojourn times of $Y^{(n,\delta)}$ and $\bar Y^{(n,\delta)}$, respectively. 
We first notice that $\bar{\mathcal{S}}_i^{(n,\delta)}\geq\mathcal{S}_i^{(n,\delta)}$ 
for every $i$.

Given $T,\eta,\delta>0$, we take $n_0,\Delta$ satisfying (\ref{S:clock}) and 
we may conclude that outside an event of probability at most $\eta$, we have that
\begin{eqnarray}\label{ineqlemma}
\max(\bar{\mathcal{S}}_i^{(n,\delta)}-\mathcal{S}_i^{(n,\delta)})\leq\!\!\!\!\sum_{x\in d_n^{-1}\mathbb{Z}\cap(\mathcal{T}_\delta^{(n)})^c}
\!\!\!\!g_n(V_{x+d_n^{-1}}-V_x)\phi_n(\Delta,x)\bar T_{xd_n}(n\Delta),
\end{eqnarray}
where the max is taken over all sojourn times of $\delta$-traps visited by $Y^{(n,\delta)}$ during $[0,T]$.

As seen before, the right side of (\ref{ineqlemma}) goes to 0 in probability by first taking $n\rightarrow\infty$ 
and after $\delta\downarrow0$. This result and the arbitrariness of $\eta$ yield the result. $\square$\\

\noindent{\bf Proof of Theorem \ref{bivagingtheorem}.} For simplicity, let us take $\lambda=\theta\mu$ 
and replace $\mu$ by $a_n$ as defined in (\ref{an}). We have that for every $n\geq1$ and $T>0$
\begin{eqnarray}\nonumber
&{\displaystyle \mathcal{\bar R}(\theta a_n,a_n)=
\int_0^\infty f(t)P(Y^{(n)}_t=Y^{(n)}_{t(1+\theta)})\,dt}&\\\label{age1}
&{\displaystyle =\int_0^T f(t) P(Y_t^{(n)}= Y_{t(1+\theta)}^{(n)})\,dt+g_n(T)}&
\end{eqnarray}
where $f$ is the probability density function of $\mathbb{T}$ and 
\begin{equation}
 g_n(T)=\int_T^\infty f(t) P(Y_t^{(n)}= Y_{t(1+\theta)}^{(n)})\,dt\leq P(\mathbb{T}>T)\label{age2}
\end{equation} 
for every $n\geq1$.

From Lemma \ref{lemma1aging}, it follows that 
\begin{eqnarray*}
 &&\lim_{n\rightarrow\infty}\int_0^T f(t) P(Y_t^{(n)}= Y_{t(1+\theta)}^{(n)})\,dt\\
&=&\lim_{\delta\downarrow0}\lim_{n\rightarrow\infty}
\int_0^T f(t) P(\bar{Y}_t^{(n,\delta)} = \bar{Y}_{t(1+\theta)}^{(n,\delta)})\,dt.
\end{eqnarray*}

By using the above result, Lemma \ref{lemma2aging} and the following weak convergence 
$\lim_{\delta\downarrow0}\lim_{n\rightarrow\infty}Y_t^{(n,\delta)}\stackrel{d}{=}Z_{S_t^{-1}}$
(under the $J_1$ metric), we get
\begin{eqnarray}\nonumber
&{\displaystyle \lim_{n\rightarrow\infty}\int_0^T f(t)P(Y^{(n)}_t=Y^{(n)}_{t(1+\theta)})\,dt}
={\displaystyle \lim_{n\rightarrow\infty}E\bigg(\int_0^T f(t)1\{Y^{(n)}_t=Y^{(n)}_{t(1+\theta)}\}\bigg)\,dt}\\
\nonumber 
&=\displaystyle{ 
E\bigg(\int_0^T f(t) 1\{Z_{S_t^{-1}}=Z_{S_{t(1+\theta)}^{-1}}\}dt\bigg)}
={\displaystyle \int_0^T f(t) P\big(Z_{S_t^{-1}}=Z_{S_{t(1+\theta)}^{-1}}\big)\,dt}&\\
\label{finaleq}
&={\displaystyle P\big(Z_{S_1^{-1}}=Z_{S_{1+\theta}^{-1}}\big)\,P(\mathbb{T}\leq T),}&
\end{eqnarray}
where the third equality follows from Proposition \ref{selfsimilar}, which in particular implies that 
$P\big(Z_{S_t^{-1}}=Z_{S_{t(1+\theta)}^{-1}}\big)$
does not depend on $t$. The result follows from~(\ref{age1}), (\ref{age2}), (\ref{finaleq}),
and the arbitrariness of $T$. $\square$

\section*{Appendix: Results on quasistable processes}

We here present some results on quasistable processes, including one that states that the function $\mathcal R(\cdot)$ is non trivial. 
Let us introduce a notation that will be used below. Let $B$ and $C$ be Borel sets on $\mathbb R$. We define $B+C=\{x+y:\, x\in B,\, y\in C\}$. 
Let $Y_t=Z_{S_t^{-1}}$, $t\geq0$. For a fixed realization of $V$, 
define $\mathcal T\equiv \bigcup_{\delta>0}\mathcal T_\delta$, the set of traps.

\begin{proposition}\label{propappendix}
For all $t>0$, we have that $\mathbb P(Y_t\in \mathcal T)=1$ for almost every $V$.
\end{proposition}

\noindent{\bf Proof.}  We first show that for any $t>0$ we have
\begin{eqnarray}\label{nontrivial1}
\int_0^t\mathbb P(Y_s\in\mathcal T)ds=t.
\end{eqnarray}

For $\delta>0$ let $Y^{(\delta)}$ be as in Section~\ref{intenv} above (see paragraph right below Remark~\ref{fin}).
Arguing similarly as in the proof of Lemma \ref{lemma1aging}, we have that, for fixed $t,\eta,\delta>0$, there exists $\Delta>0$ such that
$$\int_0^t\mathbb P(Y_s\neq Y_s^{(\delta)})ds\leq \mathbb E\bigg(\sum_{x\in\mathcal T_\delta^c}(V_x-V_{x-})\phi(\Delta,x)\bigg)+\eta t.$$

Using the fact that $\sum_{x\in\mathcal T_\delta^c}(V_x-V_{x-})\phi(\Delta,x)\stackrel{p}{\longrightarrow}0$ as $\delta\downarrow0$ and the 
arbitrariness of $\eta$, we get that
$$\lim_{\delta\downarrow0}\int_0^t\mathbb P(Y_s\neq  Y_s^{(\delta)})ds=0.$$

Hence it follows that
\begin{eqnarray*}
\int_0^t\mathbb P(Y_s\in\mathcal T)ds&=&\lim_{\delta\downarrow0}\int_0^t\mathbb P(Y_s\in\mathcal T_\delta, Y_s= Y_s^{(\delta)})ds\\
&=&\lim_{\delta\downarrow0}\int_0^t\mathbb P(Y_s= Y_s^{(\delta)})ds=t,
\end{eqnarray*}
and (\ref{nontrivial1}) is established. 

For arbitrary $d>0$, define the set $B_d=\{s\in(0,d):\, \mathbb P(Y_s\in\mathcal T)=1\}$.  
From (\ref{nontrivial1}), we have that $\mathcal L(B_d)=d$ (we recall that $\mathcal L$ is the Lebesgue measure).

It can be seen using the Markov property that if $t$ and $s$ belong to $B_d$, then $B_d + B_d\subset B_{2d}$.
Since the sum of sets of positive Lebesgue measure contains an interval (see Theorem 4.1 from \cite{leth1988}), 
we get that $B_{2d}$ contains a subinterval for all $d>0$. Let $B=\bigcup_{d>0}B_d$. We have that\\

\noindent 1) From the Markov property, $C+D\subset B$ for any subsets $C$ and $D$ from $B$;\\

\noindent 2) From the above result, $B$ contains a subinterval of $[0,d]$ for all $d>0$.\\

Let $I_1\equiv [d_1,d_2]$ be a subinterval of $[0,d]$ obtained from 2), with 
$d_1 < d_2$. 
From 1), we have that $I_1+I_1 =[2d_1,2d_2]\subset B$. So, it follows that 
$I_2\equiv I_1\cup (I_1+I_1)=[d_1,2d_2]\subset B$. Inductively, we find that
$I_{n+1}\equiv I_1\cup (I_n+I_1)=[d_1,nd_2]\subset B$ for $n\in\mathbb N$. 
It follows that $[d_1,\infty)\subset B$.
Since $d$ is arbitrary, we conclude that $B=(0,\infty)$. $\square$

\begin{lemma}\label{lemmaappendix}
For all $t>0$ and $x\in\mathcal T$, we have that $\mathbb P_x(Y_t=x)>0$ for almost every $V$, 
where the subscript $x$ in $\mathbb P_x(\cdot)$ means that $Y_0=x$.
\end{lemma}

\noindent{\bf Proof.} Let $x\in\mathcal T$ and $d>0$. We have that for all $d>0$
\begin{eqnarray}\label{nontrivial2}
\int_0^d\mathbb P_x(Y_s=x)ds&=&\mathbb E_x\bigg(\int_0^d1\{Y_s=x\}ds\bigg)\nonumber\\
&=&\mathbb E_x[\phi(x,S^{-1}_d)] (V_x-V_{x-})>0,
\end{eqnarray}
since $\mathbb P_x$-almost surely $S_t>0$ and $\phi(x,t)>0$ for all $t>0$, where $\phi(\cdot)$ is the 
local time of $Z_t$ (see first result of~\cite{stone1963}, on page 632). 

Let $F_d=\{s\in(0,d):\,\mathbb P_x(Y_s=x)>0\}$. From (\ref{nontrivial2}), we get that 
$\mathcal L(F_d)>0$ for all $d>0$. Let $F=\bigcup_{d>0}F_d$. Arguing similarly
as in the proof of Proposition \ref{propappendix}, we find that $F=(0,\infty)$, so proving
the desired result. $\square$

\begin{corollary}\label{coro}
The function $\mathcal R(\theta)=P(Z_{S^{-1}_1}=Z_{S^{-1}_{1+\theta}})=P(Y_1=Y_{1+\theta})$ is non trivial.
\end{corollary}

\noindent{\bf Proof.} From Proposition \ref{propappendix} and Lemma \ref{lemmaappendix}, we immediately obtain that
$\mathbb P(Y_1=Y_{1+\theta})=\sum_{x\in\mathcal T}\mathbb P(Y_1=x)\mathbb P_x(Y_\theta=x)>0$ for almost every $V$ 
and for all $\theta>0$. Therefore $\mathcal R(\theta)=E[\mathbb P(Y_1=Y_{1+\theta})]>0$, where the expectation is 
taken with respect to $V$. That  $\mathcal R(\theta)$ is not constant may be verified by showing that for almost every $V$ we have that 
$\mathbb P_x(Y_t=y)\longrightarrow0$ as $t\to\infty$ for every $x,y\in\mathcal T$. We leave the details as an exercise.
$\square$





\end{document}